\documentclass{article}

\usepackage{amsmath,amsfonts,amssymb,amsthm}
\usepackage{xy}
\usepackage[bookmarks=false,breaklinks=true]{hyperref}
\usepackage{float}
\usepackage{graphicx}
\usepackage[percent]{overpic}
\usepackage{color}

\definecolor{myurlcolor}{rgb}{0,0,0.7}
\usepackage{hyperref}
\hypersetup{colorlinks,
linkcolor=myurlcolor,
citecolor=myurlcolor,
urlcolor=myurlcolor}
\usepackage[mathscr]{euscript}

\input xy
\xyoption{all}
\newcommand{\maps}{\colon}    %correct symbol for colon in f: X -> Y
\newcommand{\Z}{{\mathbb Z}}  %integers
\newcommand{\Q}{{\mathbb Q}}  %rational numbers
\newcommand{\R}{{\mathbb R}}  %real numbers
\newcommand{\C}{{\mathbb C}}  %complex numbers
\renewcommand{\H}{{\mathbb H}}  %quaternions
\newcommand{\CP}{\mathbb{C}\mathrm{P}} % complex projective space

\newcommand{\SO}{{\rm SO}}      %special orthogonal group
\newcommand{\SU}{{\rm SU}}      %special unitary group
\newcommand{\Sp}{{\rm Sp}}      %symplectic group
\newcommand{\A}{\mathrm{A}}     % A series, alternating groups
\newcommand{\E}{\mathrm{E}}     % E series
\renewcommand{\S}{\mathrm{S}}   % symmetric groups

\newcommand{\I}{\mathscr{I}} % Klein's icosahedral function
\newcommand{\II}{\mathbb{I}} % icosians

\theoremstyle{definition}

        \newcommand{\be}{\begin{equation}}
        \newcommand{\ee}{\end{equation}}
        \newcommand{\ba}{\begin{eqnarray}}
        \newcommand{\ea}{\end{eqnarray}}
        \newcommand{\ban}{\begin{eqnarray*}}
        \newcommand{\ean}{\end{eqnarray*}}
        \newcommand{\barr}{\begin{array}}
        \newcommand{\earr}{\end{array}}

\title{From the Icosahedron to $\E_8$}

\author{John C.\ Baez\\[.5em]
{\small Department of Mathematics} \\[-.3em]
{\small  University of California}\\[-.3em]
{\small Riverside, California 92521, USA} \\
\small and \\
{\small Centre for Quantum Technologies}  \\[-.3em]
{\small National University of Singapore} \\[-.3em]
{\small Singapore 117543}  \\
\small  baez@math.ucr.edu 
}

\date{\small \today}

\begin{document}

\maketitle

\begin{abstract}
\noindent
The regular icosahedron is connected to many exceptional objects in mathematics.  Here we describe two constructions of the $\E_8$ lattice from the icosahedron.  One uses a subring of the quaternions called the ``icosians'', while the other uses du Val's work on the resolution of Kleinian singularities.   Together they link the golden ratio, the quaternions, the quintic equation, the 600-cell, and the Poincar\'e homology 3-sphere.  We leave it as a challenge to the reader to find the connection between these two constructions.
\end{abstract}

\section*{}

In mathematics, every sufficiently beautiful object is connected to all others.   Many exciting adventures, of various levels of difficulty, can be had by following these connections.   
Take, for example, the icosahedron---that is, the \emph{regular} icosahedron, one of the five Platonic solids.  Starting from this it is just a hop, skip and a jump to the $\E_8$ lattice, a wonderful pattern of points in 8 dimensions!  As we explore this connection we shall see that it also ties together many other remarkable entities: the golden ratio, the quaternions, the quintic equation, a highly symmetrical 4-dimensional shape called the 600-cell, and a manifold called the Poincar\'e homology 3-sphere.   

Indeed, the main problem with these adventures is knowing where to stop.  The story we shall tell is just a snippet of a longer one involving the McKay correspondence and quiver representations.  It would be easy to bring in the octonions, exceptional Lie groups, and more.  But it can be enjoyed without these esoteric digressions, so let us introduce the protagonists without further ado.

The icosahedron has a long history.  According to a comment in Euclid's \textsl{Elements} it was discovered by Plato's friend Theaetetus, a geometer who lived from roughly 415 to 369 BC.  Since Theaetetus is believed to have classified the Platonic solids, he may have found the icosahedron as part of this project.   If so, it is one of the earliest mathematical objects discovered as part of a classification theorem.   In any event, it was known to Plato: in his \textsl{Timaeus}, he argued that water comes in atoms of this shape.

The icosahedron has 20 triangular faces, 30 edges, and 12 vertices.  We can take the vertices to be the four points
\[                    (0 , \pm 1 , \pm \Phi)   \]
and all those obtained from these by cyclic permutations of the coordinates, where
\[             \Phi = \frac{\sqrt{5} + 1}{2} \]
is the golden ratio.  Thus, we can group the vertices into three orthogonal ``golden rectangles'': rectangles whose proportions are $\Phi$ to $1$.  

\begin{center}
\includegraphics[scale = 0.3]{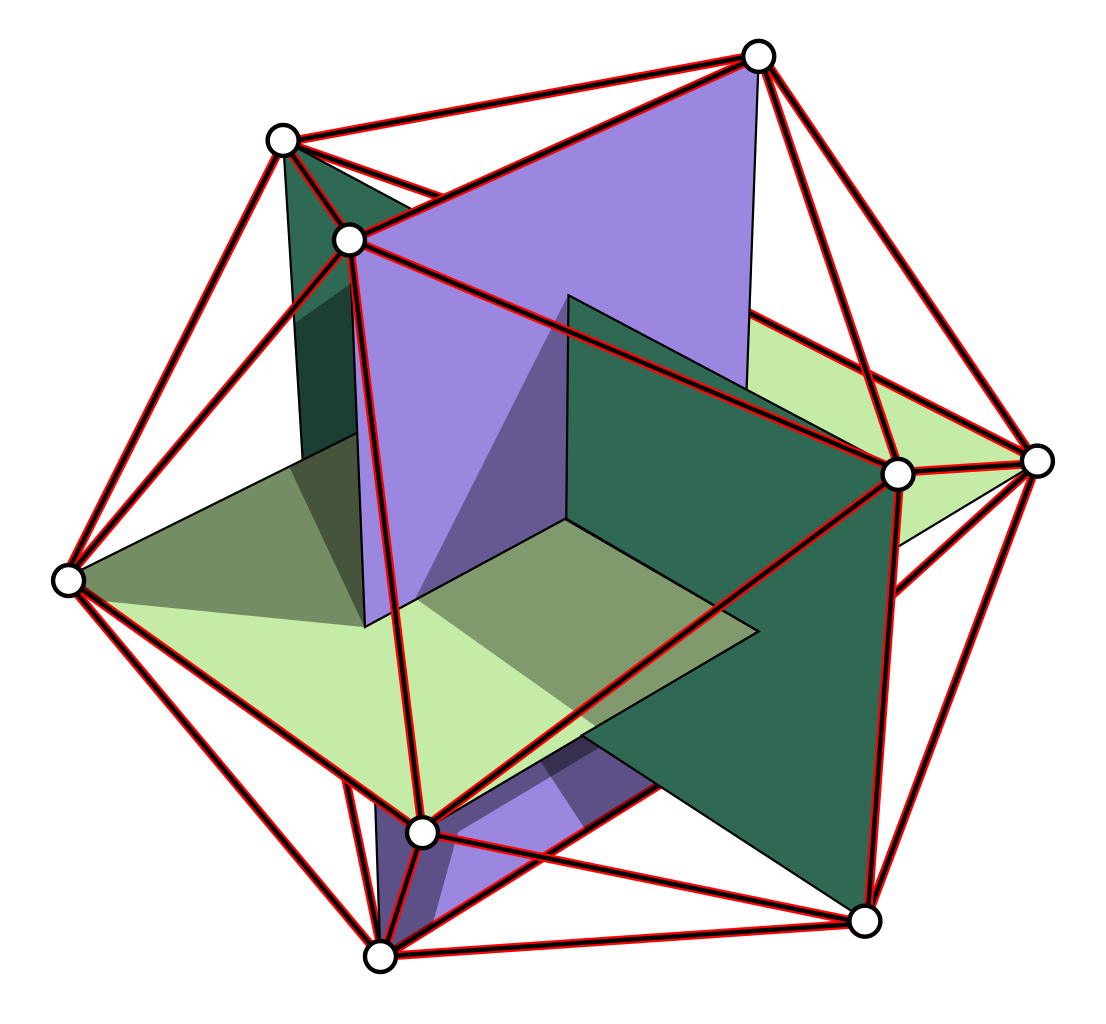}
\break
Figure 1: Icosahedron with three golden rectangles.
\end{center}

In fact, there are five ways to do this.  The rotational symmetries of the icosahedron permute these five ways, and any nontrivial rotation gives a nontrivial permutation.  The rotational symmetry group of the icosahedron is thus a subgroup of $\S_5$.   Moreover, this subgroup has $60$ elements.  After all, any rotation is determined by what it does to a chosen face of the icosahedron: it can map this face to any of the 20 faces, and it can do so in 3 ways.  The rotational symmetry group of the icosahedron is thus a 60-element subgroup of $\S_5$.  Group theory therefore tells us that it must be the alternating group $\A_5$.

The $\E_8$ lattice is harder to visualize than the icosahedron, but still easy to characterize.   Take a bunch of equal-sized spheres in 8 dimensions.  Get as many of these spheres to touch a single sphere as you possibly can. Then, get as many to touch \emph{those} spheres as you possibly can, and so on.  Unlike in 3 dimensions, where there is ``wiggle room'', you have no choice about how to proceed, except for an overall rotation and translation.  The balls will inevitably be centered at points of the $\E_8$ lattice!  

We can also characterize the $\E_8$ lattice as the one giving the densest packing of spheres among all lattices in 8 dimensions.   This packing was long suspected to be optimal even among those that do not arise from lattices---but this fact was proved only in 2016, by the young mathematician Maryna Viazovska \cite{Viazovska}.  

We can also describe the $\E_8$ lattice more explicitly.   In suitable coordinates, it consists of vectors for which:
\begin{enumerate}
\item the components are either all integers or all integers plus $\textstyle{\frac{1}{2}}$, and
\item the components sum to an even number. 
\end{enumerate}
This lattice consists of all integral linear combinations of the $8$ rows of this matrix:
\[\left( \begin{array}{rrrrrrrr}
1&-1&0&0&0&0&0&0 \\
0&1&-1&0&0&0&0&0 \\
0&0&1&-1&0&0&0&0 \\
0&0&0&1&-1&0&0&0 \\
0&0&0&0&1&-1&0&0 \\
0&0&0&0&0&1&-1&0 \\
0&0&0&0&0&1&1&0 \\
-\frac{1}{2}&-\frac{1}{2}&-\frac{1}{2}&-\frac{1}{2}&-\frac{1}{2}&-\frac{1}{2}&-\frac{1}{2}&-\frac{1}{2} 
\end{array} \right)
\]
The inner product of any row vector with itself is $2$, while the inner product of distinct row vectors is either $0$ or $-1$.   Thus, any two of these vectors lie at an angle of either $90^\circ$ or $120^\circ$.      If we draw a dot for each vector, and connect two dots by an edge when the angle between their vectors is $120^\circ$, we get this pattern:
\[ \xymatrix{ \bullet   \ar@{-}[r] & \bullet \ar@{-}[r] &\bullet \ar@{-}[r] & \bullet \ar@{-}[r] &\bullet \ar@{-}[r] \ar@{-}[d] & \bullet \ar@{-}[r] &\bullet \\
&&&& \bullet
}\]
This is called the $\E_8$ Dynkin diagram.  In the first part of our story we shall find the $\E_8$ lattice hiding in the icosahedron; in the second part, we shall find this diagram.  The two parts of this story must be related---but the relation remains mysterious, at least to this author.  

\subsection*{The Icosians}

The quickest route from the icosahedron to $\E_8$ goes through the fourth dimension.  The symmetries of the icosahedron can be described using certain quaternions; the integer linear combinations of these form a subring of the quaternions called the ``icosians'', but the icosians can be reinterpreted as a lattice in 8 dimensions, and this is the $\E_8$ lattice.   Let us see how this works, following the treatment in Conway and Sloane's \textsl{Sphere Packings, Lattices and Groups} \cite{ConwaySloane,Dechant}.
  
The quaternions, discovered by Hamilton, are a 4-dimensional algebra 
\[                \H = \{a + bi + cj + dk \colon \; a,b,c,d\in \R\}  \]
with multiplication given as follows:
\[     i^2 = j^2 = k^2 = -1, \]
\[    i j = k = - j i  \textrm{ and cyclic permutations}  .\]
It is a normed division algebra, meaning that the norm
\[    |a + bi + cj + dk| = \sqrt{a^2 + b^2 + c^2 + d^2} \]
obeys $|q q'| = |q| |q'|$ for all $q,q' \in \H$.  The unit sphere in $\H$, a 3-sphere, is therefore a group.  We shall sloppily call this group $\SU(2)$, since its elements can be identified with $2 \times 2$ unitary complex matrices with determinant $1$.  A more correct name for this group would be $\Sp(1)$.    But whatever we call it, this group acts as rotations of 3-dimensional Euclidean space, since we can see any point in $\R^3$ as a ``purely imaginary'' quaternion $x = bi + cj + dk$, and $qxq^{-1}$ is then a purely imaginary quaternion of the same norm for any $q \in \SU(2)$.   Indeed, this action gives a homomorphism
\[           \alpha \maps \SU(2) \to \SO(3) \]
where $\SO(3)$ is the group of rotations of $\R^3$.   This homomorphism is two-to-one, since $\pm 1 \in \SU(2)$ both act trivially.   So, we say $\SU(2)$ is a ``double cover'' of the rotation group $\SO(3)$.

We can thus take any Platonic solid, look at its group of rotational symmetries, get a subgroup of
$\SO(3)$, and take its double cover in $\SU(2)$.  If we do this starting with the icosahedron, we see that the $60$-element group $\A_5 \subset \SO(3)$ is covered by a 120-element group $\Gamma \subset \SU(2)$, called the ``binary icosahedral group''.   Of course the elements of $\Gamma$ lie on the unit sphere in the quaternions---but it turns out that they are the vertices of a 4-dimensional regular polytope: a 4-dimensional cousin of the Platonic solids!  This polytope looks like a  ``hypericosahedron'', but it is usually called the ``600-cell'', since it has 600 tetrahedral faces.

\begin{center}
\includegraphics[scale = 1.1]{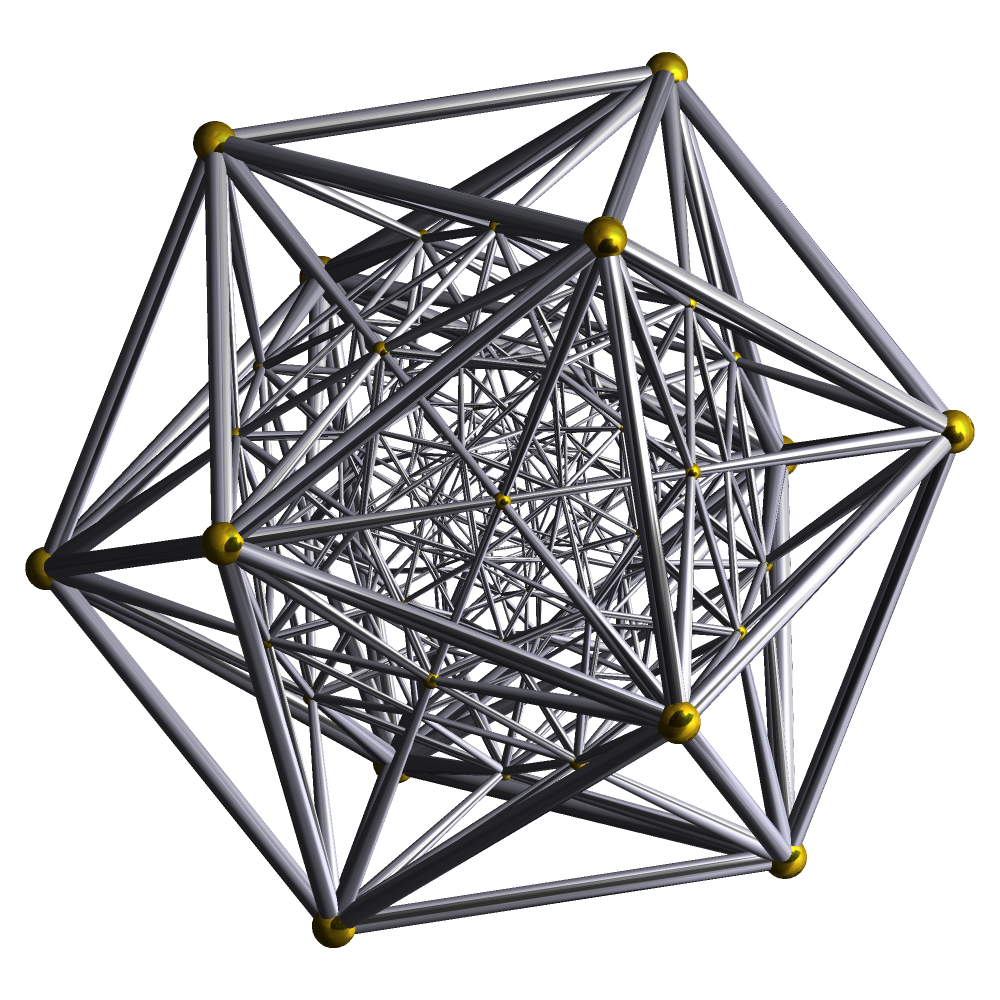}
\break Figure 2: The 600-cell projected down to 3 dimensions, \break
drawn using Robert Webb's \href{http://www.software3d.com/Stella.php}{Stella software}.
\end{center}

Explicitly, if we identify $\H$ with $\R^4$, 
the elements of $\Gamma$ are the points
\[   (\pm \textstyle{\frac{1}{2}}, \pm \textstyle{\frac{1}{2}},\pm \textstyle{\frac{1}{2}},\pm \textstyle{\frac{1}{2}}) , \quad  (\pm 1, 0, 0, 0) , \quad
    \textstyle{\frac{1}{2}} (\pm \Phi, \pm 1 , \pm \textstyle{\frac{1}{\Phi}}, 0 ),
\]
and those obtained from these by even permutations of the coordinates.
Since these points are closed under multiplication, if we take integral linear combinations of them we get a subring of the quaternions:
\[             \II = \{ \sum_{q \in \Gamma} a_q  q  : \; a_q \in \Z \}  \subset \H .\]
 Conway and Sloane \cite{ConwaySloane} call this the ring of ``icosians''.  The icosians are not a lattice in the quaternions: they are dense.  However, any icosian is of the form $a + bi + cj + dk$ where $a,b,c$, and $d$ live in the ``golden field'' 
\[            \Q(\sqrt{5}) = \{ x + \sqrt{5} y : \; x,y \in \Q\}. \]
Thus we can think of an icosian as an 8-tuple of rational numbers.  Such 8-tuples form a lattice
in 8 dimensions.

In fact we can put a norm on the icosians as follows.  For $q \in \II$ the usual quaternionic
norm has
\[            |q|^2 =  x + \sqrt{5} y \]
for some rational numbers $x$ and $y$, but we can define a new norm on $\II$ by setting
\[        \|q\|^2 = x + y .\]
With respect to this new norm, the icosians form a lattice that fits isometrically in 8-dimensional Euclidean space.  This is none other than $\E_8$!

\subsection*{Klein's Icosahedral Function}

Not only is the $\E_8$ lattice hiding in the icosahedron; so is the $\E_8$ Dynkin diagram.  The space of all regular icosahedra of arbitrary size centered at the origin has a singularity, which corresponds to a degenerate special case: the icosahedron of zero size.  If we resolve this singularity in a minimal way we get eight Riemann spheres, intersecting in a pattern described by the $\E_8$ Dynkin diagram!

This part of our story starts around 1884, with Felix Klein's \textsl{Lectures on the Icosahedron} \cite{Klein}.  In this work he inscribed an icosahedron in the Riemann sphere, $\CP^1$.  He thus got the icosahedron's symmetry group, $\A_5$, to act as conformal transformations of $\CP^1$---indeed, rotations.  He then found a rational function of one complex variable that is invariant under all these transformations.  This function equals $0$ at the centers of the icosahedron's faces, 1 at the midpoints of its edges, and $\infty$ at its vertices. 

\begin{center}
\includegraphics[scale = 0.4]{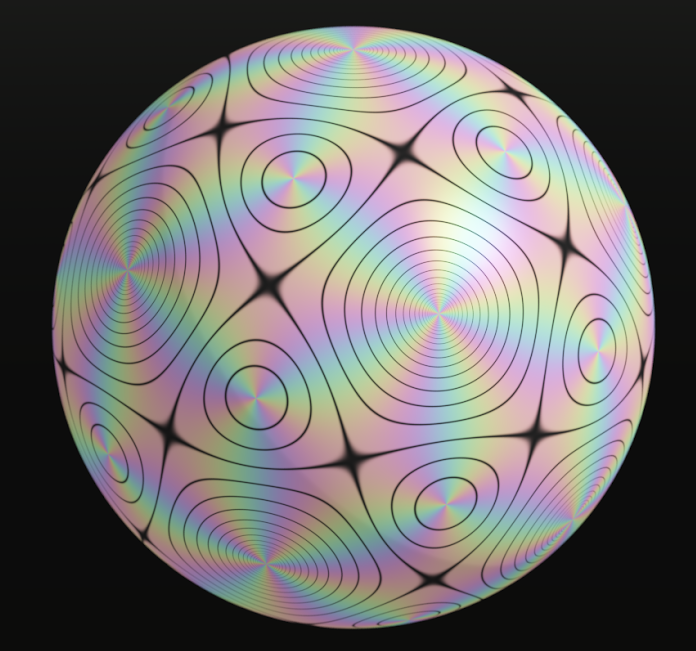} 
\break
Figure 3: Klein's icosahedral function drawn by Abdelaziz Nait Merzouk.  \\ 
The color shows its phase, while the contour lines show its magnitude.
\end{center}

We can think of Klein's icosahedral function as a branched cover of the Riemann sphere by itself with 60 sheets:
\[               \I \maps \CP^1 \to \CP^1 .\]
Indeed, $\A_5$ acts on $\CP^1$, and the quotient space $\CP^1/\A_5$ is isomorphic to $\CP^1$ again.  The function $\I$ gives an explicit formula for the quotient map $\CP^1 \to \CP^1/\A_5 \cong \CP^1$.

Klein managed to reduce solving the quintic to the problem of solving the equation $\I(z) = w$ for $z$.   A modern exposition of this result is Shurman's \textsl{Geometry of the Quintic} \cite{Shurman}.  For a more high-powered approach, see Nash \cite{Nash}.   Unfortunately, neither of these treatments avoids complicated calculations.  But our interest in Klein's icosahedral function here does not come from its connection to the quintic: instead, we want to see its connection to $\E_8$. 

For this we should actually construct Klein's icosahedral function.   To do this, recall that the Riemann sphere  $\CP^1$ is the space of 1-dimensional linear subspaces of $\C^2$.   Let us work directly with $\C^2$.  While $\SO(3)$ acts on $\CP^1$, this comes from an action of this group's double cover $\SU(2)$ on $\C^2$.    As we have seen, the rotational symmetry group of the icosahedron, $\A_5 \subset \SO(3)$, is double covered by the binary icosahedral group $\Gamma \subset \SU(2)$.    To build an $\A_5$-invariant rational function on $\CP^1$, we should thus look for $\Gamma$-invariant homogeneous polynomials on $\C^2$.

It is easy to construct three such polynomials:
\begin{itemize}
\item $V$, of degree $12$, vanishing on the 1d subspaces corresponding to icosahedron vertices.
\item $E$, of degree $30$, vanishing on the 1d subspaces corresponding to icosahedron edge midpoints.
\item $F$, of degree $20$, vanishing on the 1d subspaces corresponding to icosahedron face centers. 
\end{itemize}
Remember,  we have embedded the icosahedron in $\CP^1$.  Each point in $\CP^1$ is a 1-dimensional subspace of $\C^2$, so each icosahedron vertex determines such a subspace, and
there is a nonzero linear function on $\C^2$, unique up to a constant factor, that vanishes on this subspace. The icosahedron has $12$ vertices, so we get $12$ linear functions this way. Multiplying them gives $V$, a homogeneous polynomial of degree $12$ on $\C^2$ that vanishes on all the subspaces corresponding to icosahedron vertices.  The same trick gives $E$, which has degree $30$ because the icosahedron has $30$ edges, and $F$, which has degree $20$ because the icosahedron has $20$ faces.   

A bit of work is required to check that $V,E$ and $F$ are invariant under $\Gamma$, instead of changing by constant factors under group transformations.  Indeed, if we had copied this construction using a tetrahedron or octahedron, this would not be the case.  For details, see Shurman's free book \cite{Shurman} or van Hoboken's nice thesis \cite{Hoboken}.

Since both $F^3$ and $V^5$ have degree $60$, $F^3/V^5$ is homogeneous of degree zero, so it defines a rational function $\I \maps \CP^1 \to \CP^1$.  This function is invariant under $\A_5$ because $F$ and $V$ are invariant under $\Gamma$.  Since $F$ vanishes at face centers of the icosahedron while $V$ vanishes at vertices, $\I = F^3/V^5$ equals $0$ at face centers and $\infty$ at vertices.  Finally, thanks to its invariance property, $\I$ takes the same value at every edge center, so we can normalize $V$ or $F$ to make this value $1$.   

Thus, $\I$ has precisely the properties required of Klein's icosahedral function!  And indeed, these properties uniquely characterize that function, so that function is $\I$.

\subsection*{The Appearance of $\E_8$}

Now comes the really interesting part.    Three polynomials on a 2-dimensional space must obey 
a relation, and $V,E,$ and $F$ obey a very pretty one, at least after we normalize them correctly:
\[     V^5 + E^2 + F^3 = 0. \]  
We could guess this relation simply by noting that each term must have the same degree.  Every 
$\Gamma$-invariant polynomial on $\C^2$ is a polynomial in $V, E$ and $F$, and indeed
\[         \C^2 / \Gamma \cong  \{ (V,E,F) \in \C^3 \colon \; V^5 + E^2 + F^3 = 0 \} . \]
This complex surface is smooth except at $V = E = F = 0$, where it has a singularity.  And hiding in this singularity is $\E_8$!

To see this, we need to ``resolve'' the singularity.   Roughly, this means that we find a smooth complex surface $S$ and an onto map
\[ \xymatrix{  S \ar[d]^(.4){\pi}  \\ \C^2/\Gamma } \]
that is one-to-one away from the singularity\footnote{More precisely, if $X$ is an algebraic variety with singular points $X_{\mathrm{sing}} \subset X$, $\pi \maps S \to X$ is a ``resolution'' of $X$ if $S$ is smooth, $\pi$ is proper, $\pi^{-1}(X - X_{\textrm{sing}})$ is dense in $S$, and $\pi$ is an isomorphism between $\pi^{-1}(X - X_{\mathrm{sing}})$ and $X - X_{\mathrm{sing}}$.  For more details see \cite{Lamotke}.}.   There are many ways to do this, but one ``minimal'' resolution, meaning that all others factor uniquely through this one:
\[ \xymatrix{ S' \ar[dr]_(.4){\pi'} \ar@{-->}[r]^{\exists !} & S \ar[d]^(.4){\pi}  \\ & \C^2/\Gamma } \]

What sits above the singularity in this minimal resolution?  Eight copies of the Riemann sphere $\CP^1$, one for each dot here:
\[ \xymatrix{ \bullet   \ar@{-}[r] & \bullet \ar@{-}[r] &\bullet \ar@{-}[r] & \bullet \ar@{-}[r] &\bullet \ar@{-}[r] \ar@{-}[d] & \bullet \ar@{-}[r] &\bullet \\
&&&& \bullet
}\]
Two of these $\CP^1$s intersect in a point if their dots are connected by an edge: otherwise they are disjoint.  

This amazing fact was discovered by Patrick du Val in 1934 \cite{DuVal}. Why is it true?  Alas, there is not enough room in the margin, or even the entire page, for a full explanation.  The books by Kirillov \cite{Kirillov} and Lamotke \cite{Lamotke} fill in the details.   But here is a clue.  The $\E_8$ Dynkin diagram has ``legs'' of lengths $5, 2$ and $3$:
\[ \xymatrix{ 5  \ar@{-}[r] & 4 \ar@{-}[r] & 3 \ar@{-}[r] & 2 \ar@{-}[r] & 1 \ar@{-}[r] \ar@{-}[d] & 2 \ar@{-}[r] & 3 \\
&&&& 2
} \]
On the other hand, 
\[  \A_5 \cong \langle v, e, f \mid v^5 = e^2 = f^3 = vef = 1 \rangle \]
where in terms of the icosahedon's symmetries:
\begin{itemize}
\item    $v$ is a $1/5$ turn around some vertex of the icosahedron,
\item   $e$ is a $1/2$ turn around the center of an edge touching that vertex, 
\item   $f$ is a $1/3$ turn around the center of a face touching that vertex,
\end{itemize}
and we must choose the sense of these rotations correctly to obtain $vef = 1$.  To get a presentation of the binary icosahedral group we drop one relation:
\[ \displaystyle{  \Gamma \cong \langle v, e, f | v^5 = e^2 = f^3 = vef \rangle } .\]
The dots in the $\E_8$ Dynkin diagram correspond naturally to conjugacy classes in $\Gamma$, not counting the conjugacy class of the central element $-1$.  Each of these conjugacy classes, in turn, gives a copy of $\CP^1$ in the minimal resolution of $\C^2/\Gamma$. 

Not only the $\E_8$ Dynkin diagram, but also the $\E_8$ lattice, can be found in the minimal
resolution of $\C^2/\Gamma$.   Topologically, this space is a 4-dimensional manifold.  
Its real second homology group is an 8-dimensional vector space with an inner product given by the intersection pairing.  The integral second homology is a lattice in this vector space spanned by the 8 copies of $\CP^1$ we have just seen---and it is a copy of the $\E_8$ lattice \cite{KS}.

But let us turn to a more basic question: what is $\C^2/\Gamma$ like as a topological space?  To tackle this, first note that we can identify a pair of complex numbers with a single quaternion, and this gives a homeomorphism
\[         \C^2/\Gamma \cong \H/\Gamma \] 
where we let $\Gamma$ act by right multiplication on $\H$.  So, it suffices to understand $\H/\Gamma$.  

Next, note that sitting inside $\H/\Gamma$ are the points coming from the unit sphere in $\H$. These points form the 3-dimensional manifold $\SU(2)/\Gamma$, which is called the ``Poincar\'e homology 3-sphere'' \cite{KS}.   This is a wonderful thing in its own right:  Poincar\'e discovered it as a counterexample to his guess that any compact 3-manifold with the same homology as a 3-sphere is actually diffeomorphic to the 3-sphere, and it is deeply connected to $\E_8$.   But for our purposes, what matters is that we can think of this manifold in another way, since we have a diffeomorphism
\[    \SU(2)/\Gamma \cong \SO(3)/\A_5.  \]
The latter is just \emph{the space of all icosahedra inscribed in the unit sphere in 3d space}, where we count two as the same if they differ by a rotational symmetry. 

This is a nice description of the points of $\H/\Gamma$ coming from points in the unit sphere of $\H$.  But every quaternion lies in \emph{some} sphere centered at the origin of $\H$, of possibly zero radius.   It follows that $\C^2/\Gamma \cong \H/\Gamma$ is the space of \emph{all} icosahedra centered at the origin of 3d space---of arbitrary size, including a degenerate icosahedron of zero size.   This degenerate icosahedron is the singular point in $\C^2/\Gamma$.    This is where $\E_8$ is hiding.

Clearly much has been left unexplained in this brief account.  Most of the missing details can be found in the references.  But it remains unknown---at least to this author---how the two constructions of $\E_8$ from the icosahedron fit together in a unified picture.

Recall what we did.   First we took the binary icosahedral group $\Gamma \subset \H$, took integer linear combinations of its elements, thought of these as forming a lattice in an 8-dimensional rational vector space with a natural norm, and discovered that this lattice is a copy of the $\E_8$ lattice.  Then we took $\C^2/\Gamma \cong \H/\Gamma$, took its minimal resolution, and found that the integral 2nd homology of this space, equipped with its natural inner product, is a copy of the $\E_8$ lattice.  From the same ingredients we built the same lattice in two very different ways!  How are these constructions connected?  This puzzle deserves a nice solution.  

\section*{Acknowledgements}

I thank Tong Yang for inviting me to speak on this topic at the Annual General Meeting of the Hong Kong Mathematical Society on May 20, 2017, and Guowu Meng for hosting me at the HKUST while I prepared that talk.  I also thank the many people, too numerous to accurately list, who have helped me understand these topics over the years.

\end{document}